\documentclass[11pt]{article}
\usepackage{amsmath}
\usepackage{amsthm}
\usepackage{amsfonts}
\usepackage{setspace}
\usepackage{fullpage}
\usepackage{amssymb}
\usepackage{enumitem}

\usepackage{comment}

\usepackage{tikz}
\usepackage{float}
\usepackage{caption}

\newtheorem{theorem}{Theorem} 
\newtheorem{proposition}{Proposition} 
 
\newtheorem{definition}{Definition}
\newtheorem{problem}{Problem}

\newcommand{\cH}{\mathcal{H}}

\newcommand{\EH}{{\rm EH}}

\date{}
\newcounter{num}

\title{Homogeneous sets in hypergraphs with forbidden order-size pairs}

\author{Maria Axenovich\thanks{Karlsruhe Institute of Technology, Karlsruhe, Germany, 
	\texttt{maria.aksenovich@kit.edu}.	Research supported in part by the DFG grant FKZ AX 93/2-1.} \and Dhruv Mubayi\thanks{University of Illinois at Chicago, Chicago, USA, 	 \texttt{mubayi@uic.edu}. Research partially supported by NSF grants DMS-1763317, 1952767, 2153576, a Humboldt research award and a Simons fellowship.} \and Lea Weber\thanks{Karlsruhe Institute of Technology, Karlsruhe, Germany, \texttt{lea.weber@kit.edu}.}}

\begin{document}

\maketitle

\begin{abstract}The well-known Erd\H{o}s-Hajnal conjecture states that for any  graph $F$, there exists $\epsilon>0$ such that  every $n$-vertex graph $G$ that contains no induced copy of $F$ has a homogeneous   set of size at least $n^{\epsilon}$. We consider a variant of the Erd\H{o}s-Hajnal problem for hypergraphs where we forbid a family of hypergraphs described by their orders and sizes. 
For graphs, we observe that if we forbid induced subgraphs on $m$ vertices and $f$ edges for any positive $m$ and $0\leq f \leq \binom{m}{2}$, then we obtain large homogeneous sets. For triple systems, in the first nontrivial case $m=4$, for every $S \subseteq \{0,1,2,3,4\}$, we give bounds on the minimum size of a homogeneous set in a triple system where the number of edges spanned by every four vertices is not in $S$. For all $S$ we determine if the growth rate is polylogarithmic. Several open problems remain. 
\end{abstract}

\section{Introduction}
For an integer $r\geq 2$, an $r$-{\it graph} or $r$-uniform hypergraph is a pair $H=(V, E)$, where $V=V(H)$ is the set of vertices and  $E=E(H) \subseteq \binom{V}{r}$ is the set  of edges.  A $2$-graph is simply a graph.  A {\it homogeneous set} is a set of vertices that is either a clique or a coclique (independent set). For an $r$-graph $H$, let $h(H)$ be the size of a largest homogeneous set.  Given $r$-graphs $F, H$, say that $H$ is $F$-{\it free}  if $H$ contains no isomorphic copy of $F$ as an induced subgraph.
We say that an $r$-graph $F$ has the {\it Erd\H{o}s-Hajnal-property} or simply {\it {\rm EH}-property} if there is a constant $\epsilon=\epsilon_F>0$ such that every $n$-vertex $F$-free $r$-graph $H$ satisfies $h(H) \geq n^{\epsilon}$.  A conjecture of Erd\H{o}s and Hajnal~\cite{EH} states that  any $2$-graph has the EH-property.  The conjecture remains open, see for example  a survey by Chudnovsky~\cite{C}, as well as \cite{APS, BLT, FPS}, to name a few.   When $F$ is a fixed graph and $G$ is an $F$-free $n$-vertex graph, Erd\H{o}s and Hajnal proved that $h(G) \ge 2^{c\sqrt{\log n}}$. This was recently improved to $h(G) \ge 2^{c\sqrt{\log n\log\log n}}$ by   Buci\'{c}, Nguyen, Scott, and Seymour~\cite{BNSS}.

The Erd\H{o}s-Hajnal conjecture fails for $r$-graphs, $r\geq 3$, already when $F$  is a clique of size $r+1$.  Indeed,  well-known results on off-diagonal hypergraph Ramsey numbers show that there are $n$-vertex $r$-graphs that do not have a clique on $r+1$ vertices and do not have cocliques on $f_r(n)$ vertices, where $f_r$ is an iterated logarithmic function (see~\cite{MS} for the best known results). 
 Moreover, the following   result (Claim 1.3. in \cite{GT})   tells us exactly  which $r$-graphs, $r\geq 3$,  have the EH-property. Here $D_2$ is the unique $3$-graph on $4$ vertices with exactly $2$ edges.

\begin{theorem}[Gishboliner and Tomon~\cite{GT}]\label{GT}
Let $r\geq 3$. If $F$ is an $r$-graph on at least $r+1$ vertices and $F\neq D_2$, then there is an $F$-free  $r$-graph  $H$ on $n$ vertices such that $h(H) =(\log n)^{O(1)}$.
\end{theorem}

It is natural to consider the EH-property for families of $r$-graphs instead of a single $r$-graph. 
In this paper, we consider families determined by a given set of orders and sizes. Several special cases of this have been extensively studied over the years (see, e.g.~\cite{EH1}).
 For  $0\leq f \leq \binom{m}{r}$,   we call an $r$-graph $F$  on $m$ vertices and $f$ edges an $(m,f)$-\emph{graph} and we call the pair $(m,f)$ the {\it order-size pair} for $F$.
 Say that $H$ is $(m,f)$-free if it contains no induced copy of an $(m,f)$-graph.
 If $Q=\{(m_1, f_1), \ldots, (m_t, f_t)\}$, say that $H$ is $Q$-free if $H$ is $(m_i,f_i)$-free for all $i=1, \ldots, t$.
 
 \begin{definition}
 	Given $r \ge 2$ and 
 		 $Q=\{(m, f_1), \ldots, (m, f_t)\}$, let $h(n,Q)=h_r(n,Q)$ be the minimum of $h(H)$, taken over all  $n$-vertex $Q$-free $r$-graphs $H$.  Say that  $Q$  has the \EH-{\it property} if there exists $\epsilon=\epsilon_Q>0$ such that $h(n, Q) >n^{\epsilon}$. 
 \end{definition}
 
 For example $h_3(n, \{(4,0), (4,2)\}) = k$ means that any $n$-vertex $3$-graph in which  any $4$ vertices induce $1$, $3$, or $4$ edges has a homogenous set of size $k$, and there is an $r$-graph $H$ as above with $h(H) = k$.  We may omit the subscript $r$ in the notation $h_r(n, Q)$ if it is obvious from context. When $Q=\{(m,f)\}$ we use the simpler notation $h(n,m,f)$ instead of $h(n, \{(m,f)\})$.  Let us make two simple observations:  \begin{equation} \label{subset}
 	h_r(n, Q) \le h_r(n, Q') \qquad \text{ if} \qquad  Q\subseteq Q',
 	\end{equation}

 \begin{equation} \label{complement} h_r(n, Q) = h_r(n, \overline Q) \qquad \text{ where } \qquad 
 	\overline Q = \left\{\left(m, {m \choose r}-f\right): (m,f) \in Q\right\}.
 	\end{equation}
 Our first result concerns 2-graphs, where we show that forbidding a single order-size pair already guarantees large homogeneous sets.
 
 \begin{proposition}\label{graph}
 For any integers  $m,f$ with  $m\geq 2$ and $0\leq f\leq \binom{m}{2}$ there exists $c>0$ such that 
 $h_2(n, m, f) > c \, n^{1/(m-1)}$.
  \end{proposition}

It seems a challenging problem to give good upper bounds on $h_2(n, m, f)$. For example, determining $h_2(n,m,{m \choose 2})$ is equivalent to determining off-diagonal Ramsey numbers.  
 
 Our remaining results are in the hypergraph case $r=3$ and $m=4$. 
  We shall be considering sets $Q$ of pairs $(4,i)$ for $i\in \{0, 1, 2, 3, 4\}$. 
 We do not need to consider sets $Q$ that contain both $(4,0)$ and $(4,4)$ because Ramsey's theorem guarantees that for sufficiently large $n$ we cannot avoid both of them. 
 Using complementation~(\ref{complement}), this leaves us with the following sets: 
 \begin{itemize}
 	\item{} $\{(4,0)\}$, ~$\{(4,1)\}$, ~$\{(4,2)\}$; \\
 	\item{}  $ \{ (4,0), (4,1) \}$,~  $ \{ (4,0), (4,2) \}$,   ~$ \{ (4,0), (4,3) \}$,  ~$ \{ (4,1), (4,2) \}$, ~$ \{ (4,1), (4,3) \}$; ~~ and \\
 	\item{} $\{ (4,0), (4,1), (4,2) \}$,  ~$\{ (4,0), (4,1), (4,3) \}$, ~$\{ (4,0), (4,2), (4,3) \}$, ~$\{ (4,1), (4,2), (4,3) \}$.
 \end{itemize} 
We address $h(n, Q)$ for each of these choices of $Q$.

We quickly obtain bounds for the first case using results in Ramsey theory (note again that $h(n,4,f)=h(n,4,4-f)$). Recall that the Ramsey number  $R_k(s,t)$ is the minimum $n$ such that every red/blue edge-coloring of the complete $n$-vertex $k$-graph yields either a monochromatic red $s$-clique or a monochromatic blue $t$-clique. It is known~\cite{CFS} that $2^ {c t\log t} \leq  R_3 (4, t)  \leq   2^{c' t^2 \log t}$. This yields positive constants $c$ and $c'$, such that
	\begin{equation} \label{40lower}c' \left(\frac{\log n}{\log\log n}\right)^{1/2} < h_3(n, 4,0)  <   c\frac{\log n}{\log\log n}.
		\end{equation}

 A more recent result of Fox and He~\cite{FH} constructs $n$-vertex $3$-graphs with every four vertices spanning at most two edges and independence number at most $c\log n/\log\log n$. Together with (\ref{subset}) this yields positive a constant  $c$, such that
 	\begin{equation} \label{41upper} h_3(n, 4,1)  \le
 		h_3(n, \{(4,0), (4,1)\})
 		<  c \frac{\log n}{\log\log n}.
 \end{equation}
 
 For the remaining cases when $|Q|=1$ we obtain  bounds using recent results by Fox and He~\cite{FH} and  by Gishboliner and Tomon~\cite{GT}. Recall that $f(n) = (1+o(1))g(n)$ means that there is a function $e(n)$ such that $\lim_{n \rightarrow \infty} e(n) \rightarrow 0$ and $f(n) = g(n) + e(n)g(n)$.

\begin{proposition} \label{thm:41} There are  positive constants $c_1, c_2$  such that 
	  \begin{equation} \label{eq:41}   h_3(n, 4,1 ) > c_1 \, \left(\frac{\log n}{\log\log n}\right)^{1/2}    \end{equation}
and 	  \begin{equation} \label{eq:42} n^{c_2}< h_3(n,4,2) < (1+o(1))n^{1/2}. \end{equation}
	  	\end{proposition} 
 It is unclear if either bound for $h(n,4,1)$ above represents the correct order of magnitude, but the lower bound certainly seems far off. 
 
\begin{problem} Improve the exponent $1/2$ in the lower bound on $h_3(n,4,1)$.
\end{problem}

 Our next results address  the case when $|Q|=2$. For the first case we have constants $c,c'>0$ such that
 $$c \, \frac{\log n}{\log\log n}< h_3(n,\{(4,0), (4,1)\} < c' \frac{\log n}{\log\log n}.$$
 The lower bound follows (after applying (\ref{complement})) from an old result of Erd\H os and Hajnal~\cite{EH1}. This is the first instance of a (different) conjecture of Erd\H os and Hajnal~\cite{EH1} about the growth rate of generalized hypergraph Ramsey numbers that correspond to our setting of  $h(n, Q)$, 
 where $Q = \{ (m, f), (m, f+1), \ldots, (m, \binom{m}{r}) \}$. Recent results of the second author and Razborov~\cite{MR} on this problem determine, for each $m>r \ge 4$, the minimum $f$ such that $h_r(n, Q) < c \log ^an$ for  some $a$ and  $Q=\{(m,f), \ldots, (m, \binom{m}{r})\}$. When $r=3$, the  minimum $f$ was determined by Conlon, Fox, and  Sudakov \cite{CFS} for $m$ being a power of $3$ and for growing $m$, as well as some other values.

 For the second case when $|Q|=2$, we have $h_3(n, \{(4,0),(4,2)\}) > n^c$ as follows immediately from (\ref{subset}) and (\ref{eq:42}). However, the value of $c$ obtained from~\cite{GT} is very small (less than 0.005). We improve this below to $1/5$  and also obtain bounds for the other cases. 
 
 \begin{theorem} \label{|Q|=2} There is a positive constant $c_1$  such that for $n>5$
 	
 	\begin{equation}\label{4042} h_3(n, \{(4,0),(4,2)\}) > c_1 \, n^{1/5},  \end{equation}
 	
 	\begin{equation}\label{4043} h_3(n, \{(4,0),(4,3)\}) > c_1 \, n^{1/3}, \end{equation}

 	\begin{equation}\label{4142} c_1(n \log n)^{1/3} < h_3(n, \{(4,1),(4,2)\}) < (1+o(1)){n}^{1/2},  \mbox{  and }  \end{equation}
 	
 	\begin{equation}\label{4143} \frac{1}{2}  \log n \le  h_3(n, \{(4,1),(4,3)\}) < 4(\log n)^2. \notag \end{equation}
 	\end{theorem}
 We note the upper bound 
$$h_3(n, \{(4,0),(4,2)\}) \le h_3(n, \{(4,0),(4,1), (4,2)\}) < c \sqrt{n \log n}$$ that we will see below. Apart from this we were not able to obtain  nontrivial upper bounds in (\ref{4042}) or (\ref{4043}). Improving the bounds in (\ref{4042}), (\ref{4043}) and (\ref{4142}) seems to be an interesting open problem.
 
 \begin{problem}
 	Prove or disprove that 
 	\begin{itemize}
 			\item $h_3(n, \{(4,0),(4,2)\}) = n^{1/2+o(1)}$, 
 			
 	\item  $h_3(n, \{(4,0),(4,3)\}) = n^{1+o(1)}$,  
 	
 	\item $h_3(n, \{(4,1),(4,2)\}) = n^{1/2+o(1)}$.
 	\end{itemize}
 	\end{problem}
 
 Finally, we consider the case when $|Q|=3$. If $Q=\{(4,0), (4,1), (4,2)\}$, then a $\overline Q$-free 3-graph is a partial Steiner system (STS), and it is well known~\cite{EHR, PR, DPR} that the minimum independence number of an $n$-vertex partial STS has order of magnitude $\sqrt{n \log n}$. Thus
 $h_3(n, Q)$ has order of magnitude  $\sqrt{n \log n}$.  If $Q=\{(4,1),(4,2), (4,3)\}$, and $n \ge 4 $, then it is a simple exercise to show that any $Q$-free 4-graph on at least four vertices is a clique or coclique and therefore $h_3(n, Q)=n$ for $n \ge 4$.
 The two remaining cases are covered below.
 \begin{theorem} \label{|Q|=3} Let $n \ge 4$. Then $h_3(n, \{(4,0),(4,2), (4,3)\}) =n-1$ and
 	
 	\begin{equation}\label{013} h_3(n, \{(4,0), (4,1), (4,3)\}) =\begin{cases}
 			 \frac{n}{2}  &\text {if  $n \equiv 0$ (mod 6)} \\
 \lceil \frac{n+1}{2}\rceil & \text {if  	$n \not\equiv 0$ (mod 6).}
 \end{cases}	  \notag 
 	    \end{equation}
 	 \end{theorem}

\medskip

In Section 2 we prove Proposition~\ref{graph} and in Section 3 we prove our results for triple systems.

\section{ Graphs}\label{graphs}

In this section we prove Proposition~\ref{graph}. 
For a graph $G$, let $\omega(G)$ and $\alpha(G)$ denote the size of a largest clique and coclique, respectively.
 \medskip

\noindent
{\bf Proof of Proposition~\ref{graph}.}
We shall use induction on $m$ with basis $m=2$. In this case $f\in \{0,1\}$.  Note that $h(n, 2, 0) =h(n, 2, 1) =n = n^1 = n^{1/(m-1)}$, 
since forbidden graphs are either a non-edge or an edge. 
Consider an $(m, f)$-free graph $G$ on $n$ vertices, $m\geq 3$, and assume that the statement of the proposition holds for smaller values of $m$.
We can also assume that $G$ is not a complete graph, an empty graph, a cycle, or the complement of a cycle, since we are done in these cases.
Consider $\Delta$ and $\overline{\Delta}$, the maximum degree of $G$ and of the complement $\overline{G}$ of $G$, respectively.
Using Brooks' theorem,  the chromatic number of $G$ and  of $\overline{G}$ is at most $\Delta$ and $\overline{\Delta}$, respectively.
Thus, $\alpha(G)\geq n/\Delta$ and  $\omega(G)\geq n/\overline{\Delta}$. Therefore, we can assume that $\Delta \geq n^{(m-2)/(m-1)}$ and $\overline{\Delta}\geq n^{(m-2)/(m-1)}$, otherwise we are done. Thus, there is a vertex with at least $n^{(m-2)/(m-1)}$ edges incident to it and there is a vertex with at least $n^{(m-2)/(m-1)}$ non-edges incident to it.

Assume first that $f\leq m-1$. Consider a vertex $v$ with at least $n^{(m-2)/(m-1)}$ non-edges incident to it, i.e., with a set $X$ of vertices each non-adjacent to $v$, $|X|\geq n^{(m-2)/(m-1)}$. Since $G$ is $(m,f)$-free, $G[X]$ is $(m-1, f)$-free. Thus, by induction $h(G) \geq h(G[X]) \geq |X|^{1/(m-2)}\geq 
n^{1/(m-1)}.$

Now assume that $f\geq m$.  Consider a vertex $v$ with at least $n^{(m-2)/(m-1)}$ edges incident to it, i.e., with a set $X$ of vertices each adjacent to $v$, $|X|\geq n^{(m-2)/(m-1)}$. Since $G$ is $(m,f)$-free, $G[X]$ is $(m-1, f-(m-1))$-free.  Thus, by induction 
$h(G) \geq h(G[X]) \geq |X|^{1/(m-2)}\geq 
n^{1/(m-1)}.$ \qed

\section{Triple systems}  \label{r=3}

In this section we prove Proposition~\ref{thm:41}, Theorem~\ref{|Q|=2} and Theorem~\ref{|Q|=3}.
We will need the following notions and result for our proofs.
For an $r$-graph $H$ and  one of its vertices $v$, we define the {\it link graph} of $v$ to be the $(r-1)$-graph $L(v)$  whose vertex set is $V(H)\setminus \{v\}$ and edge set is $\{e \subseteq V(H)\setminus \{v\}: e \cup \{v\} \in E(H)\}$.
When denoting edges in $3$-graphs, we often shall  omit parentheses and commas, for example instead of writing $\{x,y,z\}$ we simply shall write $xyz$.
For a $2$-graph $G$,  let $L(G)$ be the $3$-graph with vertex set $V(G) \cup \{v\}$, $v\not \in V(G)$ and edge set $\{uvw: uw\in E(G)\}$. 
Finally, when we consider a $3$-graph $H$, the link graph of a vertex $u\in V(H)$  {\it restricted to a vertex set} $S$, denoted $L_S(u)$ is a graph on vertex set $S$ and edge set $\{vw:  v, w\in S, uvw \in E(H)\}$. A clique on $s$ vertices is denoted $K_s$.

We shall use the following theorem.

\begin{theorem}[Fox and He \cite {FH}, Thm. 1.4]\label{thm:FH}
For all $t, s\geq 3$, any $3$-graph on more than $(2t)^{st}$ vertices contains either a coclique on $t$ vertices or $L(K_s)$. 
\end{theorem}

\subsection{Forbidden sets of size $1$}

{\bf Proof of Proposition~\ref{thm:41}.}
First we consider the case $Q=\{(4,1)\}$.

To prove the lower bound on $h(n, 4, 1)$, we shall consider the complementary setting and an arbitrary $n$-vertex $(4,3)$-free $3$-graph $H$. We shall apply Theorem \ref{thm:FH} with largest possible $t=s$ such that $ (2t)^{st}<n$. 
In this case $t=s \geq c(\log n / \log \log n)^{1/2}$. If  $H$  has a coclique of size $t$, then $h(H) \geq t$ and we are done. Otherwise  $H$ contains a subgraph isomorphic to  $L=L(K_s)$.  Let $V(L)= \{v\} \cup V$, where all edges are incident to $v$, $v\not\in V$.  Note that $V$ induces a clique in $H$, otherwise $v$ and three vertices of $V$ not inducing an edge give a $(4,3)$-subgraph.  Thus, $h(H) \geq s-1$. In each case $h(H) \geq c(\log n / \log \log n)^{1/2}$.

Now, we consider the case  $Q=\{(4,2)\}$.  The lower bound on $h(n, 4, 2)$ follows from a  result of Gishboliner and Tomon~\cite{GT}. The  upper bound  is obtained by taking  an affine plane of order $q$. More precisely, given a sufficiently large $n$, choose a prime $q$ such that $n^{1/2} < q \le n^{1/2} + n^{0.29}$; such $q$ exists by density results about primes (see, e.g., ~\cite{BHP}). Let $A(2,q)$ be the affine plane of order $q$.  Let $H$ be the 3-graph whose vertex set is some $n$-element subset of the point set of $A(2,q)$, and whose edge set is the set of triples that are contained in some line in $A(2,q)$. Let $S$ be a set of four vertices in $H$. If two lines each contain at least three points in $S$, then they have two points in common, which is impossible, hence at most one line contains at least three points in $S$. This means that $S$ induces 0, 1 or 4 edges, and consequently, $H$ is $(4,2)$-free. The largest clique in $H$ is the vertex set of a line, and has size at most $q$. The largest coclique in $H$ is a cap set in $A(2,q)$ which is well known to have size at most $q+2$. Hence $h(H) \le q+2 < n^{1/2} +  n^{0.3}$ for sufficiently large $n$.
\qed

\subsection{Forbidden sets of size $2$}

We will need the following special cases of  results of de-Caen~\cite{DC} on the hypergraph Tur\'an problem and of Kostochka, Mubayi, and Verstra\"ete~\cite{KMV} on independent sets in sparse hypergraphs. 

\begin{theorem} [de-Caen~\cite{DC}] \label{thm:DC}
	Suppose that $n>k\ge 3$ and $H$ is an $n$-vertex 3-graph with more   than $(1-{k-1 \choose 2}^{-1})(n^3/6)$ edges. Then $H$ contains a clique of size $k$.
	\end{theorem}

\begin{theorem} [Kostochka, Mubayi, and Verstra\"ete~\cite{KMV}] \label{thm:KMV}
	Suppose that  $H$ is an $n$-vertex 3-graph  in which every pair of vertices lies in at most $d$ edges, where $0<d<n/(\log n)^{27}$. Then $H$ has an independent set of size at least
	$c \sqrt{(n/d)\log (n/d)}$ where $c$ is an absolute constant. 
\end{theorem}

{\bf Proof of Theorem~\ref{|Q|=2}.} 

{\it Case $Q=\{(4,0), (4,2)\}$.} \\ 
Using complementation, we consider a $\{(4,2), (4,4)\}$-free $3$-graph $H$ on $n$ vertices. Assume $n$ is sufficiently large. We shall show that $h(H) \geq Cn^{1/5}$, for some constant $C>0$. 
For  a vertex $v$ in $H$, let $K$ be a clique in the link graph $L(v)$ of $v$. Then $K$ is a coclique in $H$, for an edge within $K$ in $H$ together with $v$ yields a 4-clique in $H$. We will use this observation repeatedly.  Suppose that the complement of $H$ has $(1-\gamma)(n^3/6)$ edges for some $0<\gamma < 1$ and $k\ge 3$ is defined
via 
$$\frac{1}{{k \choose 2}} \le \gamma < \frac{1}{{k-1 \choose 2}}.$$
Then  by Theorem~\ref{thm:DC}, $H$ has a coclique of size at least $k$. If $k>n^{1/5}$, then we are done so assume from now that $k<n^{1/5}$. As $H$ has at least $\gamma n^3/6 - n^2/2$ edges, by averaging, $H$ has two vertices $v,w$ whose common neighborhood $S$ has size at least $\gamma n-4$.  If $L_S(v)$ has an induced $C_4$, then it induces a $4$-clique in $H$, for otherwise we obtain a $(4,2)$-subgraph in $H$. Hence $L_S(v)$ has no induced $C_4$ and by known results (see, e.g.~\cite{G}) it has a homogeneous set $T$ of size at least $c |S|^{1/3}$. 
If $T$ is a clique in $L_S(v)$, then by our observation, $T$ is a coclique in $H$. If $T$ is a coclique in $L_S(v)$, then  $T$ is a clique in $L_S(w)$ for otherwise
we obtain a $(4,2)$-subgraph in $H$ with $v,w$ and two vertices in $T$. Again the observation implies that $T$ is a coclique in $H$. Hence in both cases $T$ is a coclique in $H$ and 
$h(H) \ge |T| \ge c|S|^{1/3} \ge (c/2) (\gamma n)^{1/3}$. Since $k<n^{1/5}$, we have $\gamma> n^{-2/5}$ and $h(H) >(c/2)n^{1/5}$ completing the proof.  \\

{\it Case $Q=\{(4,0), (4,3)\}$.} \\ 
We shall again consider the complementary case.  Suppose that $H$ is a  $3$-graph on $n$ vertices 
that is $\{(4,1), (4,4)\}$-free. We will prove that  $h(H) \geq n^{1/3}$. Let $y$ be an arbitrary vertex of $H$ and consider the link graph  $L(y)$ of $y$.

Assume that there  is an induced  $2K_2$ in $L(y)$, i.e., that there is a set $X$ of four vertices inducing exactly two disjoint edges  in $L(y)$. 
Any three vertices in $X$ form an edge in $H$, otherwise  
these three vertices and $y$ span exactly one edge in $H$, a contradiction. 
Thus, $X$ spans exactly $4$ edges in $H$, a contradiction.  Thus, $L(y)$ is $2K_2$-free.
In the graph case it is known, that $2K_2$ has the Erd\H{o}s-Hajnal property, and in particular that any $n$-vertex  graph with no induced $2K_2$ contains a homogeneous set of size  $cn^{1/3}$ (see e.g.~\cite{G}). Thus, $h(L(y))\geq cn^{1/3}$.

Note that a $3$-clique in $L(y)$ is not an edge in $H$,  since otherwise there is a $4$-clique in $H$.  Similarly, a $3$-coclique in $L(y)$ is not an edge in $H$, since otherwise  together with $y$ it  induces a $(4,1)$-subgraph of $H$.  Thus, any set of vertices that is  a clique in $L(y)$ or an independent set in $L(y)$  is an independent set in $H$.  Thus, $h(H) \geq h(L(y)) \geq cn^{1/3}$ completing the proof.\\

{\it Case $Q=\{(4,1), (4,2)\}$.} \\ 
We now prove
$c_1(n \log n)^{1/3} \le h(n, \{ (4,1), (4,2) \} ) = h(n, \{ (4,3), (4,2) \} )\leq n^{1/2}+cn^{0.3}.$
The upper bound follows immediately from the construction used in the upper bound in (\ref{eq:42}) so we turn to the lower bound.
Using (\ref{complement}), consider an $n$-vertex $3$-graph $H$ that is $\{(4,2), (4,3)\}$-free where $n$ is sufficiently large. Let $u,v$ be a pair of vertices in $H$ whose common neighborhood $S$ has maximum size $d>0$. Given vertices $x,y\in S$, the edges $xyu$ and $xyv$ are both in $H$ else $\{u,v,x,y\}$ induces a $(4,2)$ or $(4,3)$-graph. Next, any three vertices $x,y,z \in S$, 
must form an edge of $H$ otherwise $\{u,x,y,z\}$ induces a $(4,3)$-graph. Therefore $S$ induces a clique in $H$ of size $d$. If $d>n^{0.4}$, say, then we are done as $h(H)\ge d$. Recalling that $n$ is large enough, we may assume that $d \le n^{0.4} < n/(\log n)^{27}$. Now Theorem~\ref{thm:KMV} yields a coclique in $H$ of size at least $c \sqrt{(n/d)\log n}$ for some positive constant $c$. Consequently, there is a constant $c'$ such that 
	$$h(H) \ge  \max \, \{d, \, c \sqrt{(n/d)\log n} \} > c' \,(n \log n)^{1/3}.$$
		Replacing $c'$ by a possibly smaller constant $c_1$ yields the result for all $n>4$.\\
		
Note that the set of maximal cliques in any $\{(4,2), (4,3)\}$-free 3-graph $H$ forms  a linear (maybe non-uniform) hypergraph $\cH$. Thus, determining $h(H)$ amounts to finding $\max\{t, |X|\}$, where $t$ is the size of a larges hyperedge and $X$ is a largest set of vertices in $\cH$ with no three in the same hyperedge.\\

{\it Case $Q=\{(4,1), (4,3)\}$.} \\ 
Finally, we prove
$\frac{1}{2} \log n \leq h(n, \{ (4,1), (4,3) \} ) \leq 4 (\log n)^2$. For the lower bound  let $H$ be an $n$-vertex $Q$-free  $3$-graph. Pick a vertex $v$ in $H$ and consider its link graph $L(v)$. Since $R_2(t,t) < 4^{t-1}$ (see Erd\H os and Szekeres~\cite{ES}),  we see that $L(v)$ has a clique or coclique  $K$ of size at least $\frac{1}{2} \log n$.  In the first case,  $K$ is a clique in $H$, else we find a $(4,3)$-subgraph in $H$, and in the second case,  $K$ is a coclique in $H$, else we find a $(4,1)$-subgraph in $H$.
 
We now turn to the upper bound.
Let $\chi$ be a red/blue coloring of an $n$-vertex complete graph on vertex set $V$ in which every monochromatic clique has size at most $2\log n$.  Such a coloring exists by the classical result of Erd\H os~\cite{E47}.  Let $H$ be the $3$-graph  on vertex set $V$ whose edge set consists of all triples of vertices that induce a triangle with one or  three red edges under $\chi$.

Consider four vertices $x, y, z,$ and $w$ of $H$ and assume that $xyz$ is an edge in $H$. Then the triangle $xyz$ has one or three red edges under $\chi$. Assume that $xy$ is red. 
We need to treat two cases when  $xz$ and $yz$ are blue and when  $xz$ and $yz$ are red.
In each of these cases, consider the fourth vertex $w$ and possible colors on the edges from $w$ to $x, y$ and $z$. In each of these cases $\{x,y,z,w\}$ induces exactly two or exactly four edges. 
Thus, any four vertices of $H$ induce none, two, or four edges. So, $H$ is $Q$-free.
Consider a homogeneous set $S$ in $H$.  If it is a clique, all triangles with vertices in $S$ have exactly one or three red edges under $\chi$.  Thus, the graph induced by $S$ is a pairwise vertex disjoint union of red cliques.  Either one of these red cliques has size at least $\sqrt{|S|}$ or, by taking a single vertex from each of these red cliques we see that there is a blue clique of size at least $\sqrt{|S|}$ under $\chi$. Since each monochromatic clique in the coloring $\chi$ has size at most $2\log n$, we have that $|S| \leq 4 (\log n)^2$.  Similarly, if $S$ is an independent set in $H$, all triangles with vertices in $S$ have exactly one or three blue edges under $\chi$ and again we get that $|S| \leq 4 (\log n)^2$.\qed

\subsection{Forbidden sets of size $3$}

 We will need the following structural characterization of $Q$-free $3$-graphs for 
$Q= \{ (4,1), (4,3), (4,4) \}$.
\begin{theorem}[Frankl and F\"uredi~\cite{FF}]\label{lem:ff_characterisation}
	Let $H$ be an $\{ (4,1), (4,3), (4,4) \})$-free $3$-graph. Then $H$ is isomorphic to one of the following $3$-graphs: 
	\begin{enumerate}
		\item A blow-up of the $6$ vertex $3$-graph $H'$ with vertex set $V(H') = [6]$ and edge set $E(H') = \{123, 124, 345, 346, 561, 562, 135, 146, 236, 245\}$. Here for the blow-up we replace every vertex of $H'$ by an independent set, and whenever we have $3$ vertices from three distinct of those sets, they induce an edge if and only if the corresponding vertices in $H'$ do.
		\item The $3$-graph whose vertices are the points of a regular $n$-gon where $3$ vertices span an edge if and only if the corresponding points span a triangle whose interior contains the center of the $n$-gon.
	\end{enumerate} 
\end{theorem}

\noindent
{\bf Proof of Theorem~\ref{|Q|=3}.}

{\it Case $Q = \{ (4,1), (4,3), (4,4) \}$.}\\
 We are to prove that
$$ h(n, \{(4,0), (4,1), (4,3)\})= h(n, Q) =\begin{cases}
	\frac{n}{2}  &\text {if $n \equiv 0$ (mod 6)} \\
	\lceil \frac{n+1}{2}\rceil & \text {if 	$n \not\equiv 0$ (mod 6)}.
\end{cases}	 $$
First, let us prove that 
the second 3-graph $H$ in Theorem~\ref{lem:ff_characterisation} has independence number exactly $\lceil{(n+1)/2)}\rceil$. Assume  the vertex set is $[n]$ and the vertices are labeled by consecutive integers in clockwise orientation. The lower bound is by taking $\lceil{(n+1)/2)}\rceil$ consecutive vertices on the $n$-gon and noting that no three of them contain the center in their interior.   For the upper bound, let us see how many  vertices can lie in an independent set containing $1$. When $n$ is odd,  the triangle formed by  $\{1, i, (n-1)/2+i\}$  contains the center and hence is  an edge. Therefore we may pair the elements of $[n]\setminus\{1\}$ as $(2, (n+3)/2), (3, (n+5)/2), \ldots, ((n+1)/2, n)$ and note that each pair can have at most one vertex in an independent set containing 1. Hence the maximum size of an independent set containing 1 is at most $(n+1)/2$ and by vertex transitivity of $H$, the independence number of $H$ is at most $(n+1)/2$. For $n$ even we consider the $n/2-1$ pairs $(2, n/2+1), (3, n/2+2), \ldots, (n/2, n-1)$ and add the vertex $n$ to get an upper bound $n/2+1=\lceil{(n+1)/2)}\rceil$.

Next we observe that the $6$-vertex 3-graph $H'$ in Theorem~\ref{lem:ff_characterisation} has independence number exactly $3$ (we omit the short case analysis needed for the proof). Hence if we blow-up each vertex of $H'$ into sets of the same size, then we obtain $n$-vertex $3$-graphs with independence number exactly $n/2$ whenever  $n \equiv 0$ (mod 6). This concludes the proof of the upper bound.

For the lower bound, let  $H$ be $Q$-free. Then by Theorem \ref{lem:ff_characterisation}, $H$ is isomorphic to one of the two graphs described in  Theorem~\ref{lem:ff_characterisation}. If $H$ is isomorphic to the second graph, then we have already shown that its independence number is at least $(n+1)/2$, so assume that 
$H$ is isomorphic to the blow-up of the $6$-vertex $10$-edge $3$-graph $H'$.   There are $10$ non-edges in $H'$. Let $V_1, \ldots, V_6$ be the blown up vertex sets.   Since every vertex $i \in [6]$ in $H'$ is contained in exactly $5$ non-edges, we obtain 
$$ 5n = 5\sum\limits_{i \in [6]} |V_i| = \sum\limits_{j_1j_2j_3 \not\in E(H)} |V_{j_1}| + |V_{j_2}| + |V_{j_3}| .$$
By the pigeonhole principle, there is a non-edge $i_1i_2i_3$, such that $ |V_{i_1}| + |V_{i_2}| + |V_{i_3}|\geq n/2$.  Our bound follows by observing that  for any  non-edge $i_1i_2i_3$ in the original $3$-graph $H'$ the set $V_{i_1} \cup V_{i_2} \cup V_{i_3}$ is an independent set. This gives an independent set of size at least $n/2$, and if 	$n \not\equiv 0$ (mod 6), then equality cannot hold throughout (a short case analysis, which we omit,  is needed to prove this)  and we obtain an independent set of size strictly greater than $n/2$ as required.\\

{\it Case $Q = \{ (4,0), (4,2), (4,3) \}$.}\\
We now prove $h(n, \{ (4,0), (4,2), (4,3) \} )  = n-1,$ for $n\geq 4$.
Let $H$ be a $3$-graph  that is a  clique on $n-1$ vertices and a single isolated vertex, then $H$ is $Q$-free, giving us the upper bound. 

 For the lower bound, let $H$ be a $Q$-free $3$-graph on $n$ vertices, $n\geq 4$. 
 Assume that $H$ is not a clique and not a coclique. We shall show that  $H$ is a clique and a single isolated vertex.  Consider a maximal clique $S$ in $H$. Since  $|S|<n$, there is  a vertex $v\in V(H)\setminus S$. From the  maximality of $S$,  $L_S(v) $ is not a clique. If $L_S(v)$ contains an edge, then we have that for some vertices $x, y, y'$, 
 $xy\in E(L_S(v))$ and $xy'\not\in E(L_S(v))$. But then $\{v, x, y, y'\}$ induces a $(4,2)$ or a $(4,3)$-graph, a contradiction. Thus, $L_S(v)$ is an empty graph, i.e., there is no edge in $H$ containing $v$ and two vertices of $S$. 
Now assume there exists a second vertex $v' \in V(H)\setminus(S \cup \{v\})$. Then by the same argument as above, $v'$ is also not contained in any edge with two vertices from $S$. 
Consider triples $vv'x$, $x\in S$. Since $|S|\ge 3$, by the pigeonhole principle there are  two vertices $x,x'\in S$ such that either
$vv'x, vv'x'\in E(H)$ or $vv'x, vv'x'\not\in E(H)$. Then $\{v,v', x, x'\}$ induces $2$ or $0$ edges respectively, a contradiction. Thus, $|S|=n-1$ and $v$ is an isolated vertex.
 \qed

\section{Concluding Remarks}

Fix integers $m>r$. Say that a set $Q$ of order size pairs $\{(m, f_1), \ldots, (m,f_t)\}$  is Erd\H os-Hajnal (EH) if there exists $\epsilon=\epsilon_Q$ such that $h_r(n, Q)>n^{\epsilon}$.  As $|Q|$ grows, the collection of $Q$-free $r$-graphs is more restrictive, and hence $h_r(n, Q)$ grows (assuming that large $Q$-free $r$-graphs are not forbidden to exist by Ramsey's theorem). 
The case when $h_r(n, Q) = \Omega(n)$ was treated by the first author and Balogh~\cite{AB} when $r=2$.
A natural question then is to ask what is the smallest $t$ such that every $Q$ of size $t$ is EH. Call this minimum value $EH_r(m)$.  Our results for $r=3$ show that for $m=4$, all $Q$ of size 3 are EH, but there are $Q$ of size 2 which are not EH. Consequently, $EH_3(4) = 3$.  

In order to further study $EH_r(m)$, we need another definition. Given integers $m\ge r\ge 3$, let $g_r(m)$ be the number of edges in an $r$-graph on $m$ vertices obtained by first taking a partition of the $m$ vertices into almost equal parts, then taking all edges that intersect each part, and then recursing this construction within each part. For example, $g_3(7)= 13$ since we start with a complete 3-partite 3-graph with part sizes $2,2,3$ and then add one edge within the part of size 3. It is known (see, e.g.~\cite{MR}) that as $r$ grows we have
$$g_r(m) = (1+o(1))\frac{r!}{r^r-r} {m \choose r}.$$
Note that $\frac{r!}{r^r-r}$ approaches 0 as $r$ grows. The second author and Razborov~\cite{MR} proved that for all fixed $m>r> 3$, there are $n$-vertex $r$-graphs
which are $Q$-free, $Q=\{(m, i): g_r(m)<i\le {m \choose r}\}$, with $h(G)= O(\log n)$. In other words, there exists $Q$ of size ${m \choose r} - g_r(m)$ which is not EH.
This proves that $EH_r(m)  \ge {m \choose r} - g_r(m)+1$. 

 Erd\H os and Hajnal~\cite{EH1} proved that for all $m > r \ge 3$, the set $Q=\{(m, i): g_r(m) \le i \le {m \choose r}\}$ is EH. In other words, they proved that every $n$-vertex $r$-graph in which every set of $m$ vertices spans less then $g_r(m)$ edges has an independent set of size at least $n^{\epsilon}$, where $\epsilon$ depends only on $r$ and $m$.  This is a particular set $Q$ of size ${m \choose r} - g_r(m) +1$
that is EH and we speculate that every other set $Q$ of this size is also EH.

\begin{problem} Prove or disprove that for all $m>r>2$,
	$$EH_r(m) = {m \choose r} - g_r(m) +1.$$
	\end{problem}
We end by noting that  $EH_3(4) = 3 = {4 \choose 3} - g_3(4) +1$.


\begin{thebibliography}{99}

\bibitem{APS} N. Alon, J.  Pach,  and  J.  Solymosi.
{\it Ramsey-type theorems with forbidden subgraphs,}
Paul Erd\H{o}s and his mathematics,
Combinatorica 21 (2001), no. 2, 155--170.

\bibitem{AB} M.  Axenovich and J. Balogh. {\it Graphs having small number of sizes on induced k-subgraphs,}  SIAM J. Discrete Math. 21, no. 1, (2007) 264--272.

\bibitem{BHP} R. C. Baker, G. Harman, and J. Pintz. {\it The difference between consecutive primes, II}. Proceedings of the London Mathematical Society 83 (2001)  (3): 532–562.

\bibitem{BLT}
N. Bousquet,  A. Lagoutte, and S. Thomass\'e. {\it The Erd\H{o}s-Hajnal conjecture for paths and antipaths}, Journal of Combinatorial Theory, Series B 113 (2015) 261--264.

\bibitem{DPR} M.  de Brandes,  K. Phelps, and V. R\"odl. {\it Coloring Steiner triple systems,} SIAM Journal on Algebraic Discrete Methods 3.2 (1982) 241--249.

\bibitem{BNSS} M. Buci\'c, T. Nguyen, A. Scott, and P. Seymour. {\it A loglog step towards Erd\H{o}s-Hajnal,} arXiv preprint arXiv:2301.10147 (2023).

\bibitem{C} M.  Chudnovsky. {\it The Erd\H{o}s-Hajnal conjecture -- a survey},  Journal of Graph Theory 75 (2014) 178--190.

\bibitem{CFS}  D. Conlon,  J. Fox,  and  B. Sudakov. {\it  Hypergraph Ramsey numbers,} J. Amer. Math. Soc. 23 (2010), no. 1, 247--266.


\bibitem{DC} D. de Caen.
{\it Extension of a theorem of Moon and Moser on complete subgraphs,}
Ars Combin. 16 (1983), 5--10.

\bibitem{E47} P. Erd\H{o}s, {\it Some remarks on the theory of graphs}, Bull. Amer. Math. Soc.53(1947), 292--294.

\bibitem{EH1} P. Erd\H{o}s and A. Hajnal. {\it On Ramsey like theorems, Problems and results, } Combinatorics (Proc. Conf. Combinatorial Math., Math. Inst., Oxford, 1972) , pp. 123--140, Inst. Math. Appl., Southendon-Sea, 1972.

\bibitem{EH} P. Erd\H{o}s and A. Hajnal. {\it Ramsey-type theorems}, Discrete Applied Mathematics 25 (1989) 37--52.

\bibitem{EHR} P. Erd\H{o}s, A. Hajnal, and B. Rothschild. {\it On chromatic number of graphs and set-systems}, Acta Math. Hungar 16 (1966) 61--99.

\bibitem{ES} P. Erd\H os and G. Szekeres. {\it A combinatorial problem in geometry}, Compositio Mathematica, 2,  (1935) 463--470.

\bibitem{FH} J. Fox and X.  He. {\it  Independent sets in hypergraphs with a forbidden link,}  
Proceedings of the London Mathematical Society. 123(4),  (2021)  384--409.

\bibitem{FPS}  J. Fox, J. Pach,  and A. Suk. {\it Erd\H{o}s-Hajnal conjecture for graphs with bounded VC-dimension, } 33rd International Symposium on Computational Geometry, Art. No. 43, 15 pp., LIPIcs. Leibniz Int. Proc. Inform., 77, Schloss Dagstuhl. Leibniz-Zent. Inform., Wadern (2017).

\bibitem{FF} P. Frankl and Z. F\"uredi. {\it An exact result for 3-graphs},  Discrete Math. 50  no. 2-3, (1984) 323--328.


\bibitem{GT} L. Gishboliner and I. Tomon. {\it On 3-graphs with no four vertices spanning exactly two edges}, arXiv:2109.04944 (2022).

\bibitem{G} A. Gy\'arf\'as. {\it Reflections on a problem of Erd\H{o}s and Hajnal,} The Mathematics of Paul Erd\H{o}s II. Springer, New York, NY, (2013) 135--141.

\bibitem{KMV} A. Kostochka, D. Mubayi,  and J. Verstra\"ete. {\it On independent sets in hypergraphs}, Random Structures and Algorithms 44 (2014), no. 2, 224.239.

\bibitem{MR} D. Mubayi and A. Razborov. {\it Polynomial to exponential transition in Ramsey theory,} Proceedings of the London Mathematical Society 122.1 (2021) 69--92.

\bibitem{MS} D. Mubayi and  A. Suk.  {\it New lower bounds for hypergraph Ramsey numbers,} Bull. Lond. Math. Soc. 50 (2018), no. 2, 189--201.

\bibitem{PR} K. Phelps and  V. R\"odl. {\it Steiner triple systems with minimum independence number,} Ars Combin 21 (1986) 167--172.

\end{thebibliography}
\end{document}